\documentclass[reqno]{amsart}%
\usepackage{amssymb}
\usepackage[colorlinks=false,hypertex]{hyperref}%,pagebackref
\date{This manuscript was commenced on 9 March 2008 and completed on 10 December 2008}
\date{}
\allowdisplaybreaks[4]
\theoremstyle{plain}
\newtheorem{thm}{Theorem}
\newtheorem{lem}{Lemma}

\theoremstyle{remark}
\newtheorem{rem}{Remark}
\DeclareMathOperator{\td}{d\mspace{-1mu}}

\begin{document}

\title[Monotonicity and convexity relating to volume of unit ball]
{Monotonicity and logarithmic convexity relating to the volume of the unit ball}

\author[F. Qi]{Feng Qi}
\address[F. Qi]{Department of Mathematics, College of Science, Tianjin Polytechnic University, Tianjin City, 300160, China}
\email{\href{mailto: F. Qi <qifeng618@gmail.com>}{qifeng618@gmail.com}, \href{mailto: F. Qi <qifeng618@hotmail.com>}{qifeng618@hotmail.com}, \href{mailto: F. Qi <qifeng618@qq.com>}{qifeng618@qq.com}}
\urladdr{\url{http://qifeng618.spaces.live.com}}

\author[B.-N. Guo]{Bai-Ni Guo}
\address[B.-N. Guo]{School of Mathematics and Informatics, Henan Polytechnic University, Jiaozuo City, Henan Province, 454010, China}
\email{\href{mailto: B.-N. Guo <bai.ni.guo@gmail.com>}{bai.ni.guo@gmail.com}, \href{mailto: B.-N. Guo <bai.ni.guo@hotmail.com>}{bai.ni.guo@hotmail.com}}

\begin{abstract}
Let $\Omega_n$ stand for the volume of the unit ball in $\mathbb{R}^n$ for $n\in\mathbb{N}$. In the present paper, we prove that the sequence $\Omega_{n}^{1/(n\ln n)}$ is logarithmically convex and that the sequence $\frac{\Omega_{n}^{1/(n\ln n)}}{\Omega_{n+1}^{1/[(n+1)\ln(n+1)]}}$ is strictly decreasing for $n\ge2$. In addition, some monotonic and concave properties of several functions relating to $\Omega_{n}$ are extended and generalized.
\end{abstract}



\thanks{The first author was partially supported by the China Scholarship Council and the Science Foundation of Tianjin Polytechnic University}

\thanks{This paper was typeset using \AmS-\LaTeX}

\maketitle

\section{Introduction}

In~\cite[Lemma~2.39]{AVV-Expos-89}, the following result was obtained: The function
\begin{equation}
\frac1x\ln\Gamma\biggl(1+\frac{x}2\biggr)
\end{equation}
is strictly increasing from $[2,\infty)$ onto $[0,\infty)$ and
\begin{equation}\label{avv-89-funct}
\lim_{x\to\infty}\biggl[\frac1{x\ln x}\ln\Gamma\biggl(1+\frac{x}2\biggr)\biggr]=\frac12,
\end{equation}
where
\begin{equation}\label{gamma-dfn}
\Gamma(x)=\int^\infty_0t^{x-1} e^{-t}\td t
\end{equation}
for $x>0$ denotes the classical Euler gamma function $\Gamma(x)$.
From this, the following conclusions were deduced in~\cite[Lemma~2.40]{AVV-Expos-89}: The sequence $\Omega_n^{1/n}$ decreases strictly to $0$ as $n\to\infty$, the series $\sum_{n=2}^\infty\Omega_n^{1/\ln n}$ is convergent, and
\begin{equation}\label{n-to-infty-Omega-n}
\lim_{n\to\infty}\Omega_n^{1/(n\ln n)}=e^{-1/2},
\end{equation}
where
\begin{equation}
\Omega_n=\frac{\pi^{n/2}}{\Gamma(1+n/2)}
\end{equation}
stands for the $n$-dimensional volume of the unit ball $\mathbb{B}^n$ in $\mathbb{R}^n$. Further, it was conjectured in~\cite[Remark~2.41]{AVV-Expos-89} that the function in~\eqref{avv-89-funct} is strictly increasing from $[2,\infty)$ onto $\bigl[0,\frac12\bigr)$ and this would imply that the sequence $\Omega_n^{1/(n\ln n)}$ is strictly decreasing for $n\ge2$.
\par
In~\cite[Theorem~1.5]{Anderson-Qiu-Proc-1997}, it was proved that the function
\begin{equation}\label{aderson-qiu-funct}
\frac{\ln\Gamma(x+1)}{x\ln x}
\end{equation}
is strictly increasing from $(1,\infty)$ onto $(1-\gamma,1)$, where $\gamma$ is the Euler-Mascheroni constant. From this, above-mentioned conjecture in~\cite[Remark~2.41]{AVV-Expos-89} were resolved in~\cite[Corollary~3.1]{Anderson-Qiu-Proc-1997}. In addition, it was conjectured in~\cite[Conjecture~3.3]{Anderson-Qiu-Proc-1997} that the function~\eqref{aderson-qiu-funct} is concave on $(1,\infty)$.
\par
In~\cite{chen-qi-oct-04-123} and \cite[Theorem~4]{X.Li.Sci.Magna-08}, the function~\eqref{aderson-qiu-funct} was proved to be strictly increasing on $(0,\infty)$.
\par
In~\cite[Section~3]{elber-laforgia-pams}, the function~\eqref{aderson-qiu-funct} was proved to be concave for $x>1$.
\par
In~\cite[Theorem~1.1]{berg-pedersen-jcam}, it was proved that
\begin{equation}\label{n-derivative-x+1}
(-1)^{n-1}\biggl[\frac{\ln\Gamma(x+1)}{x\ln x}\biggr]^{(n)}>0
\end{equation}
for $x>0$ and $n\in\mathbb{N}$. More strongly, the reciprocal of the function~\eqref{aderson-qiu-funct} was proved in~\cite[Theorem~1.4]{berg-pedersen-jcam} to be a Stieltjes transform for $x\in\mathbb{C}\setminus(-\infty,0]$, where $\mathbb{C}$ is the set of all complex numbers. Furthermore, among other things, it was directly shown in~\cite[Theorem~1.1]{pickrock} that the function~\eqref{aderson-qiu-funct} for $x\in\mathbb{C}\setminus(-\infty,0]$ is a Pick function.
\par
In~\cite[Lemma~4]{alzer-ball-ii}, it was demonstrated that the function
\begin{equation}\label{F(x)-dfn-Alzer}
F(x)=\frac{\ln\Gamma(x+1)}{x\ln(2x)}
\end{equation}
is strictly increasing on $[1,\infty)$ and strictly concave on $[46,\infty)$.
With the help of this, the double inequality
\begin{equation}\label{Theorem-2-alzer-ball-ii}
\exp\biggl(\frac{a}{n(\ln n)^2}\biggr)
\le\frac{\Omega_{n}^{1/(n\ln n)}}{\Omega_{n+1}^{1/[(n+1)\ln(n+1)]}}
<\exp\biggl(\frac{b}{n(\ln n)^2}\biggr)
\end{equation}
was turned out in~\cite[Theorem~2]{alzer-ball-ii} to be valid for $n\ge2$ if and only if
\begin{equation}
a\le\ln2\ln\pi-\frac{2(\ln2)^2\ln(4\pi/3)}{3\ln3}=0.3\dotsm \quad\text{and}\quad b\ge\frac{1+\ln(2\pi)}2=1.4\dotsm.
\end{equation}
\par
It is clear that the function in~\eqref{avv-89-funct} is equivalent to~\eqref{F(x)-dfn-Alzer}. Therefore, the above conjecture posed in~\cite[Remark~2.41]{AVV-Expos-89} was verified once again.
\par
The first aim of this paper is to extend the ranges of $x$ such that the function $F(x)$ is both strictly increasing and strictly concave respectively as follows.

\begin{thm}\label{F(x)-thm-Qi}
On the interval $\bigl(0,\frac12\bigr)$, the function $F(x)$ defined by~\eqref{F(x)-dfn-Alzer} is strictly increasing; On the interval $\bigl(\frac12,\infty\bigr)$, it is both strictly increasing and strictly concave.
\end{thm}

The second aim of this paper is, with the aid of Theorem~\ref{F(x)-thm-Qi}, to generalize the decreasing monotonicity of the sequence $\Omega_n^{1/(n\ln n)}$ for $n\ge2$, the conjecture in~\cite[Remark~2.41]{AVV-Expos-89}, and the main result in~\cite[Corollary~3.1]{Anderson-Qiu-Proc-1997}, to the logarithmic convexity.

\begin{thm}\label{unit-ball-qi-thm}
The sequence $\Omega_{n}^{1/(n\ln n)}$ is strictly logarithmically convex for $n\ge2$. Consequently, the sequence
\begin{equation}
\frac{\Omega_{n}^{1/(n\ln n)}}{\Omega_{n+1}^{1/[(n+1)\ln(n+1)]}}
\end{equation}
is strictly decreasing for $n\ge2$.
\end{thm}

In~\cite[Lemma~3]{alzer-ball-ii}, the double inequality
\begin{equation}\label{alzer-unit-ball-log}
  \frac23<\biggl[1-\frac{\ln x}{\ln(x+1)}\biggr]x\ln x\triangleq G(x)<1
\end{equation}
was verified to be true for $x\ge3$. The right-hand side inequality in~\eqref{alzer-unit-ball-log} was also utilized in the proof of the inequality~\eqref{Theorem-2-alzer-ball-ii}.
\par
The third aim of this paper is to extend and generalize the inequality~\eqref{alzer-unit-ball-log} to a monotonicity result as follows.

\begin{thm}\label{unit-ball-thm1}
The function $G(x)$ defined in the inequality~\eqref{alzer-unit-ball-log} is strictly increasing on $(0,\infty)$ with
\begin{equation}\label{2-limits}
  \lim_{x\to0^+}G(x)=-\infty\quad\text{and}\quad  \lim_{x\to\infty}G(x)=1.
\end{equation}
\end{thm}

\section{Remarks}

Before proving our theorems, we are about to give some remarks on them and the volume of the unit ball in $\mathbb{R}^n$.

\begin{rem}
It is obvious that Theorem~\ref{F(x)-thm-Qi} extends or generalizes the conjecture posed in~\cite[Remark~2.41]{AVV-Expos-89} and the corresponding conclusions obtained in~\cite[Corollary~3.1]{Anderson-Qiu-Proc-1997} and~\cite[Lemma~4]{alzer-ball-ii} respectively.
\end{rem}

\begin{rem}
For $a>1$ and $x>0$ with $x\ne\frac1a$, let
\begin{equation}\label{F(x)-dfn-alpha}
F_a(x)=\frac{\ln\Gamma(x+1)}{x\ln(ax)}.
\end{equation}
It is very natural to assert that the function $F_a(x)$ is strictly increasing on $\bigl(0,\frac1a\bigr)$ and it is both strictly increasing and strictly concave on $\bigl(\frac1a,\infty\bigr)$. More strongly, we conjecture that
\begin{equation}
(-1)^{n-1}[F_a(x)]^{(n)}>0
\end{equation}
for $x>\frac1a$ and $n\in\mathbb{N}$.
\end{rem}

\begin{rem}
From Theorem~\ref{unit-ball-thm1}, it is easy to see that the right-hand side inequality in~\eqref{alzer-unit-ball-log} is sharp, but the left-hand side inequality can be sharpened by replacing the constant $\frac23=0.666\dotsm$ by a larger number $\frac{3(2\ln2-\ln3)\ln3}{2\ln2}=0.683\dotsc$.
\end{rem}

\begin{rem}
We conjecture that
\begin{equation}
(-1)^{k-1}[G(x)]^{(k)}>0
\end{equation}
for $k\in\mathbb{N}$ on $(0,\infty)$, that is, the function $1-G(x)$ is completely monotonic on $(0,\infty)$.
\end{rem}

\begin{rem}
In~\cite{minus-one-rgmia} and its revised version~\cite{minus-one-JKMS.tex}, the reciprocal of the function $[\Gamma(x+1)]^{1/x}$ was proved to be logarithmically completely monotonic on $(-1,\infty)$. Consequently, the function
\begin{equation}
Q(x)=\biggl[\frac{\pi^{x/2}}{\Gamma(1+x/2)}\biggr]^{1/x}=\frac{\sqrt\pi\,}{[\Gamma(1+x/2)]^{1/x}}
\end{equation}
is also logarithmically completely monotonic on $(-2,\infty)$. In particular, the function $Q(x)$ is both strictly decreasing and strictly logarithmically convex on $(-2,\infty)$. Because $Q(n)=\Omega_n^{1/n}$ for $n\in\mathbb{N}$, the sequence $\Omega_n^{1/n}$ is strictly decreasing and strictly logarithmically convex for $n\in\mathbb{N}$. This generalizes one of the results in~\cite[Lemma~2.40]{AVV-Expos-89} mentioned above.
\par
Furthermore, from the logarithmically complete monotonicity of $Q(x)$, it is easy to obtain that the sequence
\begin{equation}\label{Omega-n-1-n-ratio}
\frac{\Omega_n^{1/n}}{\Omega_{n+1}^{1/(n+1)}}
\end{equation}
is also strictly decreasing and strictly logarithmically convex for $n\in\mathbb{N}$. As a direct consequence of the decreasing monotonicity of the sequence~\eqref{Omega-n-1-n-ratio}, the following double inequality may be derived:
\begin{equation}\label{Omega-n+1-n-(n+1)}
\Omega_{n+1}^{n/(n+1)}<\Omega_n\le\biggl(\frac2{\sqrt{\pi}\,}\biggr)^n\Omega_{n+1}^{n/(n+1)}, \quad n\in\mathbb{N}.
\end{equation}
When $1\le n\le4$, the right-hand side inequality in~\eqref{Omega-n+1-n-(n+1)} is better than the corresponding one in
\begin{equation}\label{Theorem-1-ball-volume-rn}
\frac2{\sqrt{\pi}\,}\Omega_{n+1}^{n/(n+1)}\le\Omega_n <\sqrt{e}\,\Omega_{n+1}^{n/(n+1)}, \quad n\in\mathbb{N}
\end{equation}
obtained in~\cite[Theorem~1]{ball-volume-rn}.
\end{rem}

\begin{rem}
In~\cite{Open-TJM-2003-Ineq-Ext.tex}, the inequality
\begin{equation}\label{yaming-ineq}
\frac{[\Gamma(x+y+1)/\Gamma(y+1)]^{1/x}}{[\Gamma(x+y+2)/\Gamma(y+1)]^{1/(x+1)}} <\sqrt{\frac{x+y}{x+y+1}}\,
\end{equation}
was confirmed to be valid if and only if $x+y>y+1>0$ and to be reversed if and only if $0<x+y<y+1$. Taking $y=0$ and $x=\frac{n}2$ in~\eqref{yaming-ineq} leads to
\begin{equation}\label{frac[Gamma(n/2+1)]}
\frac{[\Gamma(n/2+1)]^{1/n}}{[\Gamma((n+2)/2+1)]^{1/(n+2)}} =\frac{\Omega_{n+2}^{1/(n+2)}}{\Omega_n^{1/n}} <\sqrt[4]{\frac{n}{n+2}}\,,\quad n>2.
\end{equation}
Similarly, if letting $y=1$ and $x=\frac{n+1}2>1$ in~\eqref{yaming-ineq}, then
\begin{equation}
\frac{\Omega_{n+5}^{1/(n+3)}}{\Omega_{n+3}^{1/(n+1)}} <\frac1{\pi^{2/(n+1)(n+3)}}\sqrt[4]{\frac{n+3}{n+5}}\,,\quad n\ge2.
\end{equation}
\end{rem}

\begin{rem}
In~\cite{Open-TJM-2003.tex}, the following double inequality was discovered: For $t>0$ and $y>-1$, the inequality
\begin{equation}\label{open-TJM-2003-ineq}
\biggl(\frac{x+y+1}{x+y+t+1}\biggr)^a <\frac{[\Gamma(x+y+1)/\Gamma(y+1)]^{1/x}}
{[\Gamma(x+y+t+1)/\Gamma(y+1)]^{1/(x+t)}} <\biggl(\frac{x+y+1}{x+y+t+1}\biggr)^b
\end{equation}
holds with respect to $x\in(-y-1,\infty)$ if $a\ge\max\bigl\{1,\frac1{y+1}\bigr\}$ and $b\le\min\bigl\{1,\frac1{2(y+1)}\bigr\}$. Letting $t=1$, $y=0$ and $x=\frac{n}2$ for $n\in\mathbb{N}$ in~\eqref{open-TJM-2003-ineq} reveals that
\begin{equation*}
\frac{n+2}{n+4}<\frac{[\Gamma(n/2+1)]^{2/n}}
{[\Gamma((n+2)/2+1)]^{2/(n+2)}} <\sqrt{\frac{n+2}{n+4}}\,
\end{equation*}
which is equivalent to
\begin{equation}\label{sqrt-frac-n+2n+4}
\sqrt{\frac{n+2}{n+4}}\,<\frac{\Omega_{n+2}^{1/(n+2)}}{\Omega_n^{1/n}} <\sqrt[4]{\frac{n+2}{n+4}}\,,\quad n\in\mathbb{N}.
\end{equation}
When $n\ge3$, the inequality~\eqref{frac[Gamma(n/2+1)]} is better than the right-hand side inequality in~\eqref{sqrt-frac-n+2n+4}.
\par
If taking $t=1$, $y=1$ and $x=\frac{n}2-1$, then
\begin{equation}
\frac1{\pi^{2/[(n-2)n]}}\sqrt{\frac{n+2}{n+4}}\, <\frac{\Omega_{n+2}^{1/n}}{\Omega_n^{1/(n-2)}}
<\frac1{\pi^{2/[(n-2)n]}}\sqrt[8]{\frac{n+2}{n+4}}\,,\quad n\in\mathbb{N}.
\end{equation}
\par
Amazingly, replacing $t$ by $\frac12$, $y$ by $0$, and $x$ by $\frac{n}2$ in~\eqref{open-TJM-2003-ineq} results in
\begin{equation}\label{ratio-Omega-n-n+1}
\sqrt{\frac{n+2}{n+3}}\, <\frac{\Omega_{n+1}^{1/(n+1)}}
{\Omega_{n}^{1/n}} <\sqrt[4]{\frac{n+2}{n+3}}\,
\end{equation}
for $n\ge-1$. When $n>2$, this refines the inequality~\eqref{Theorem-1-ball-volume-rn} in~\cite[Theorem~1]{ball-volume-rn}.
\par
Similarly, by setting different values of $x$, $y$ and $t$ in inequalities~\eqref{yaming-ineq} and~\eqref{open-TJM-2003-ineq}, more similar inequalities as above may be derived immediately.
\end{rem}

\begin{rem}
Now it is very clear that the inequality~\eqref{Theorem-1-ball-volume-rn} obtained in~\cite[Theorem~1]{ball-volume-rn} was thoroughly strengthened by~\eqref{Omega-n+1-n-(n+1)} and~\eqref{ratio-Omega-n-n+1} together.
\end{rem}

\begin{rem}
The inequality~\eqref{ratio-Omega-n-n+1} and other related ones derived above motivate us to ask the following question: What are the best positive constants $a\ge3$, $b\le3$, $\lambda\le1$, $\mu\ge1$, $\alpha\ge2$, and $\beta\le4$ such that the inequality
\begin{equation}\label{unit-ball-open}
\sqrt[\alpha]{1-\frac\lambda{n+a}}\, <\frac{\Omega_{n+1}^{1/(n+1)}}
{\Omega_{n}^{1/n}} <\sqrt[\beta]{1-\frac\mu{n+b}}\,
\end{equation}
holds true for $n\in\mathbb{N}$?
\end{rem}

\section{Lemmas}

In order to prove our theorems, we need the following lemmas.

\begin{lem}[\cite{subadditive-qi.tex, theta-new-proof.tex-BKMS, subadditive-qi-guo-jcam.tex}]\label{comp-thm-1}
For $x\in(0,\infty)$ and $k\in\mathbb{N}$, we have
\begin{equation}\label{qi-psi-ineq-1}
\ln x-\frac1x<\psi(x)<\ln x-\frac1{2x}
\end{equation}
and
\begin{equation}\label{qi-psi-ineq}
\frac{(k-1)!}{x^k}+\frac{k!}{2x^{k+1}}< (-1)^{k+1}\psi^{(k)}(x)<\frac{(k-1)!}{x^k}+\frac{k!}{x^{k+1}}.
\end{equation}
\end{lem}

\begin{lem}[{\cite[Theorem~2]{AJMAA-063-03}}]\label{Theorem-2-AJMAA-063-03}
For $x>1$, we have
\begin{equation}\label{ajmaa-063-03-thm1}
\left(1+\frac{1}{x}\right)^{x}>\frac{x+1}{[\Gamma(x+1)]^{1/x}}.
\end{equation}
The inequality~\eqref{ajmaa-063-03-thm1} is reversed for $0<x<1$.
\end{lem}

\begin{lem}[\cite{miqschur, Mon-Two-Seq-AMEN.tex}]\label{Mon-Two-Seq-AMEN.tex-ineq}
If $t>0$, then
\begin{gather}
\frac{2t}{2+t}<\ln(1+t)<\frac{t(2+t)}{2(1+t)}; \label{log-ineq-qi}
\end{gather}
If $-1<t<0$, the inequality~\eqref{log-ineq-qi} is reversed.
\end{lem}

\section{Proofs of theorems}

\begin{proof}[Proof of Theorem~\ref{F(x)-thm-Qi}]
The inequality~\eqref{n-derivative-x+1} means that the function~\eqref{aderson-qiu-funct} is positive and increasing on $(0,\infty)$. It is apparent that the function
\begin{equation*}
  \frac{\ln x}{\ln(2x)}=\frac1{1+\ln2/\ln x}
\end{equation*}
is positive and strictly increasing on $\bigl(0,\frac12\bigr)$ and $(1,\infty)$. Therefore, the function
\begin{equation}
  F(x)=\frac{\ln\Gamma(x+1)}{x\ln x}\cdot\frac{\ln x}{\ln(2x)}
\end{equation}
is strictly increasing on $\bigl(0,\frac12\bigr)$ and $(1,\infty)$.
\par
The inequality~\eqref{ajmaa-063-03-thm1} for $0<x<1$ in Lemma~\ref{Theorem-2-AJMAA-063-03} can be rewritten as
\begin{equation}\label{ajmaa-063-03-thm1-rew}
\frac{\ln\Gamma(x+1)}{x}<(1-x)\ln(x+1)+x\ln x.
\end{equation}
\par
Direct calculation yields
\begin{equation*}
F'(x)=\frac{x\ln(2x)\psi(x+1)-[\ln(2x)+1]\ln\Gamma(x+1)}{x^2[\ln(2x)]^2} \triangleq\frac{\theta(x)}{x^2[\ln(2x)]^2}
\end{equation*}
and
\begin{equation*}
\theta'(x)=x \ln(2x)\psi'(x+1)-\frac{\ln\Gamma(x+1)}{x}.
\end{equation*}
Utilizing the left-hand side inequality in~\eqref{qi-psi-ineq} for $k=1$, the inequality~\eqref{ajmaa-063-03-thm1-rew} and Lemma~\ref{Mon-Two-Seq-AMEN.tex-ineq} leads to
\begin{align*}
  \theta'(x)&>x\biggl[\frac{1}{x+1}+\frac{1}{2(x+1)^2}\biggr]\ln(2x)-(1-x)\ln(x+1)-x\ln x\\
  &>x\biggl[\frac{1}{x+1}+\frac{1}{2(x+1)^2}\biggr]\frac{2(2x-1)}{1+2x}
  -\frac{x(1-x)(2+x)}{2(1+x)}-x\ln x\\
  &=x\biggl[\frac{2 x^4+5 x^3+8 x^2+3 x-8}{2 (x+1)^2 (2 x+1)}-\ln x\biggr]\\
  &>x\biggl[\frac{2 x^4+5 x^3+8 x^2+3 x-8}{2 (x+1)^2 (2 x+1)}-\frac{2(x-1)}{1+x}\biggr]\\
  &=\frac{x\bigl(2x^4-3x^3+4x^2+11x-4\bigr)}{2(x+1)^2(2 x+1)}\\
  &=\frac{x[2x^2(x-1)^2+x(x+1)^2+10x-4]}{2(x+1)^2(2 x+1)}\\
  &>0
\end{align*}
for $\frac12<x<1$. Hence, the function $\theta(x)$ is strictly increasing on $\bigl(\frac12,1\bigr)$. Since
$$
\theta\biggl(\frac12\biggr)=-\ln\Gamma\biggl(\frac32\biggr)=-\biggl(\frac12\ln\pi-\ln2\biggr)>0,
$$
the function $\theta(x)$, and so the function $F'(x)$, is positive, and then the function $F(x)$ is increasing on $\bigl(\frac12,1\bigr)$. In a word, the function $F(x)$ is strictly increasing on $\bigl(0,\frac12\bigr)$ and $\bigl(\frac12,\infty\bigr)$.
\par
Ready computation gives
\begin{multline*}
\frac{x^3[\ln(2x)]^3F''(x)}{2[\ln(2x)]^2+3\ln(2x)+2} =\ln\Gamma(x+1)\\
+\frac{x^2[\ln(2x)]^2\psi'(x+1)-2x\ln(2x)[\ln(2x)+1]\psi(x+1)}{2[\ln(2x)]^2+3\ln(2x)+2}
\end{multline*}
and
\begin{equation*}
\frac{\td}{\td x}\biggl\{\frac{x^3[\ln(2x)]^3F''(x)}{2[\ln(2x)]^2+3\ln(2x)+2}\biggr\} =\frac{[\ln(2x)]^2\phi(x)}{\bigl\{2[\ln(2x)]^2+3\ln(2x)+2\bigr\}^2}
\end{equation*}
for $x>\frac12$, where, by making use of the right-hand side inequality in~\eqref{qi-psi-ineq-1} and the inequality~\eqref{qi-psi-ineq} for $k=1,2$,
\begin{align*}
\phi(x)&=[2\ln(2x)+5]\psi(x+1)\\
&\quad+x\bigl\{x\bigl[2[\ln(2x)]^2+3\ln(2x)+2\bigr]\psi''(x+1) -[4\ln(2x)+3]\psi'(x+1)\bigr\}\\
&<[2\ln(2x)+5]\biggl[\ln(x+1)-\frac1{2(x+1)}\biggr]\\
&\quad+x\biggl\{-x\bigl[2[\ln(2x)]^2+3\ln(2x)+2\bigr]\biggl[\frac1{(x+1)^2}+\frac1{(x+1)^3}\biggr]\\
&\quad-[4\ln(2x)+3]\biggl[\frac1{x+1}+\frac1{2(x+1)^2}\biggr]\biggr\}\\
&=[2\ln(2x)+5]\biggl[\ln(x+1)-\frac{1}{2(x+1)}\biggr]\\
&\quad-\frac{x\bigl\{4x(x+2)[\ln(2x)]^2+2\bigl(7x^2+16x+6\bigr)\ln(2x)+10x^2+23x +9\bigr\}}{2(x+1)^3}\\
&\triangleq\varphi(x)
\end{align*}
and
\begin{align*}
-x(x+1)^4\varphi'(x)&=2x^4+10x^3+19x^2+6x+1+2(x+4)x^2[\ln(2x)]^2\\
&\quad+x\bigl(2x^3+10x^2+20x+3\bigr)\ln(2x)-2(x+1)^4\ln(x+1)\\
&\triangleq h(x),\\
h'(x)&=8x^3+34x^2+52x+7+2x(3x+8)[\ln(2x)]^2\\
&\quad+\bigl(8x^3+34x^2+56x+3\bigr)\ln(2x)-8(x+1)^3\ln(x+1),\\
xh''(x)&=24x^3+86x^2+100x+3+4x(3x+4)[\ln(2x)]^2\\
&\quad+8x\bigl(3x^2+10x+11\bigr)\ln(2x)-24x(x+1)^2\ln(x+1)\\
&\triangleq q(x),\\
q'(x)&=4\bigl\{18x^2+57x+47+(6x+4)[\ln(2x)]^2\\
&\quad+2\bigl(9x^2+23x+15\bigr)\ln(2x)-6\bigl(3x^2+4x+1\bigr)\ln(x+1)\bigr\},\\
xq''(x)&=4\bigl\{36x^2+97x+30+6[\ln(2x)]^2x-12(3x+2)\ln(x+1)x\\
&\quad+\bigl(36x^2+58x+8\bigr)\ln(2x)\bigr\},\\
x^2(x+1)q^{(3)}(x)&=8\bigl\{18x^3+53x^2+18x-11-18(x+1)x^2\ln(x+1)\\
&\quad+2\bigl(9x^3+12x^2+x-2\bigr)\ln(2x)\bigr\}\\
&\triangleq 8p(x),\\
xp'(x)&=2\bigl[27x^3+65x^2+10x-2-9(3x+2)x^2\ln(x+1)\\
&\quad+\bigl(27x^2+24x+1\bigr)x\ln(2x)\bigr],\\
(x+1)x^2p''(x)&=2\bigl[54x^4+152x^3+90x^2+3x+2+6\bigl(9x^2+13x\\
&\quad+4\bigr)x^2\ln(2x)-18\bigl(3x^2+4x+1\bigr)x^2\ln(x+1)\bigr],\\
(x+1)^2x^3p^{(3)}(x)&=2\bigl[54x^5+168x^4+146x^3+18x^2-9x-4\\
&\quad+54(x+1)^2x^3\ln(2x)-54(x+1)^2x^3\ln(x+1)\bigr],\\
(x+1)^3x^4p^{(4)}(x)&=-4\bigl(3x^5+8x^4-9x^2-19x-6\bigr)\\
&\quad\triangleq-4\lambda(x),\\
\lambda'(x)&=15 x^4+32 x^3-18 x-19,\\
\lambda''(x)&=60 x^3+96 x^2-18.
\end{align*}
Since $\lambda''(x)$ is strictly increasing for $x>\frac12$ and $\lambda''\bigr(\frac12\bigr)=\frac{27}2$, the function $\lambda''(x)>0$, so $\lambda'(x)$ is strictly increasing for $x>\frac12$. From $\lambda'\bigl(\frac12\bigr)=-\frac{369}{16}$ and $\lim_{x\to\infty}\lambda'(x)=\infty$, it follows that the function $\lambda'(x)$ has a unique zero which is a minimum of $\lambda(x)$ for $x>\frac12$. Since $\lambda\bigl(\frac12\bigr)=-\frac{549}{32}$ and $\lim_{x\to\infty}\lambda(x)=\infty$, the functions $\lambda(x)$ and $p^{(4)}(x)$ has a unique zero which is a maximum of $p^{(3)}(x)$. Due to $p^{(3)}\bigl(\frac12\bigr)=188-108\ln\bigl(\frac{3}{2}\bigr)=144.20\dotsm$ and $\lim_{x\to\infty}p^{(3)}(x)=108(1+\ln2)$, the function $p^{(3)}(x)$ is positive, thus the function $p''(x)$ is strictly increasing for $x>\frac12$. Owing to $p''\bigl(\frac12\bigr)=258-90 \ln\bigl(\frac{3}{2}\bigr)=221.50\dotsc$, the function $p''(x)$ is positive and $p'(x)$ is strictly increasing for $x>\frac12$. By virtue of $p'\bigl(\frac12\bigr)=\frac{1}{2} \bigl[181-63 \ln \bigl(\frac{3}{2}\bigr)\bigr]=77.72\dotsc$, it is easy to see that $p'(x)>0$ and $p(x)$ is strictly increasing for $x>\frac12$. Owing to $p\bigl(\frac12\bigr)=\frac{27}{4} \bigl[2-\ln\bigl(\frac{3}{2}\bigr)\bigr]=10.76\dotsc$, it is deduced that $p(x)>0$ and $q^{(3)}(x)>0$ for $x>\frac12$, consequently the function $q''(x)$ is strictly increasing for $x>\frac12$. On account of $q''\bigl(\frac12\bigr)=700-168 \ln\bigl(\frac{3}{2}\bigr)=631.88\dotsc$, we obtain that $q''(x)>0$ and $q'(x)$ is strictly increasing for $x>\frac12$. In virtue of $q'\bigl(\frac12\bigr)=320-90 \ln\bigl(\frac{3}{2}\bigr)$, we have $q'(x)>0$, and then the function $q(x)$ is strictly increasing for $x>\frac12$. Because of $q\bigl(\frac12\bigr)=\frac{155}{2}-27\ln\bigl(\frac{3}{2}\bigr)=66.55\dotsc$, we acquire that $q(x)>0$ and $h''(x)>0$, so $h'(x)$ is strictly increasing for $x>\frac12$. By $h'\bigl(\frac12\bigr)=\frac{85}{2}-27 \ln\bigl(\frac{3}{2}\bigr)=33.55\dotsc$, it is derived that $h'(x)>0$ and $h(x)$ is strictly increasing for $x>\frac12$. Since $h\bigl(\frac12\bigr)=\frac{81}{8} \bigl[1-\ln\bigl(\frac{3}{2}\bigr)\bigr]=6.01\dotsc$, we gain that $h(x)>0$ and $\varphi'(x)<0$, accordingly $\varphi(x)$ is strictly decreasing for $x>\frac12$. Due to $\varphi\bigl(\frac12\bigr)=-\frac{91}{27} +5\ln\bigl(\frac{3}{2}\bigr)=-1.34\dotsc$, it is procured that $\varphi(x)<0$, i.e., $\phi(x)<0$, hence
$$
\frac{\td}{\td x}\biggl\{\frac{x^3[\ln(2x)]^3F''(x)}{2[\ln(2x)]^2+3\ln(2x)+2}\biggr\}<0
$$
and the function
$$
\frac{x^3[\ln(2x)]^3F''(x)}{2[\ln(2x)]^2+3\ln(2x)+2}
$$
is strictly decreasing for $x>\frac12$. By
$$
\biggl\{\frac{x^3[\ln(2x)]^3F''(x)}{2[\ln(2x)]^2+3\ln(2x)+2}\biggr\}\bigg\vert_{x=1/2} =\ln\frac{\sqrt{\pi}\,}2=-0.12\dotsc,
$$
it follows that
$$
\frac{x^3[\ln(2x)]^3F''(x)}{2[\ln(2x)]^2+3\ln(2x)+2}<0,
$$
which is equivalent to $F''(x)<0$ for $x>\frac12$. As a result, the function $F(x)$ is concave on $\bigl(\frac12,\infty\bigr)$. The proof of Theorem~\ref{F(x)-thm-Qi} is complete.
\end{proof}

\begin{proof}[Proof of Theorem~\ref{unit-ball-qi-thm}]
Let
\begin{equation}
  f(x)=\biggl[\frac{\pi^{x/2}}{\Gamma(1+x/2)}\biggr]^{1/(x\ln x)}
\end{equation}
for $x>0$. Taking the logarithm of $f(x)$ gives
\begin{equation*}
\ln f(x)=\frac{\ln\pi}{2\ln x}-\frac{\ln\Gamma(1+x/2)}{x\ln x}=\frac{\ln\pi}{2\ln x} -2F\biggl(\frac{x}2\biggr),
\end{equation*}
where $F(x)$ is defined by~\eqref{F(x)-dfn-Alzer} for $x>0$ and $x\ne\frac12$. Differentiating twice gives
\begin{equation*}
[\ln f(x)]''=\frac{(\ln x+2)\ln\pi}{2x^2(\ln x)^3} -\frac12F''\biggl(\frac{x}2\biggr).
\end{equation*}
Since $F(x)$ is strictly concave for $x>\frac12$, then the function $[\ln f(x)]''$ is positive for $x>1$. As a result, the function $f(x)$ is strictly logarithmically convex for $x>1$, and so the sequence $f(n)=\Omega_{n}^{1/(n\ln n)}$ for $n\ge2$ is strictly logarithmically convex.
\par
Since the function $f(x)$ is strictly logarithmically convex for $x>1$, then $[\ln f(x)]'$ is strictly increasing for $x>1$, therefore
\begin{equation*}
\biggl[\ln\frac{f(x)}{f(x+1)}\biggr]'=[\ln f(x)-\ln f(x+1)]'=[\ln f(x)]'-[\ln f(x+1)]'<0
\end{equation*}
for $x>1$. This implies that the function $\frac{f(x)}{f(x+1)}$ is strictly decreasing for $x>1$, hence the sequence $\frac{f(n)}{f(n+1)}$ is also strictly decreasing for $n\ge2$. The proof of Theorem~\ref{unit-ball-qi-thm} is thus completed.
\end{proof}

\begin{proof}[Proof of Theorem~\ref{unit-ball-thm1}]
It is easy to see that the function $G(x)$ in~\eqref{alzer-unit-ball-log} may be rearranged as
\begin{equation}\label{f-rew}
G(x)=\frac{\ln x}{\ln(x+1)}\ln\biggl(1+\frac1x\biggr)^x.
\end{equation}
It is common knowledge that the function $\bigl(1+\frac1x\bigr)^x$ is strictly increasing and greater than $1$ on $(0,\infty)$ and tends to $e$ as $x\to\infty$. Furthermore, for $x>0$,
\begin{equation*}
\biggl[\frac{\ln x}{\ln(x+1)}\biggr]'=\frac{(x+1)\ln(x+1)-x\ln x}{x(x+1)[\ln(x+1)]^2}>0.
\end{equation*}
Thus, the function $\frac{\ln x}{\ln(x+1)}$ is strictly increasing on $(0,\infty)$ and positive on $(1,\infty)$ and, by L'H\^ospital's rule, tends to $1$ as $x\to\infty$. Therefore, the second limit in~\eqref{2-limits} is valid and the function $G(x)$ is strictly increasing on $(1,\infty)$.
The first limit in~\eqref{2-limits} can be calculated by L'H\^ospital's rule as follows
\begin{align*}
\lim_{x\to0^+}G(x)&=\lim_{x\to0^+}\frac{1-{\ln x}/{\ln(x+1)}}{1/x}\lim_{x\to0^+}\ln x\\
&=\lim_{x\to0^+} \frac{x[(x+1)\ln(x+1)-x\ln x]}{(x+1)[\ln(x+1)]^2}\lim_{x\to0^+}\ln x\\
&=\lim_{x\to0^+} \biggl\{\frac{x}{\ln(x+1)}\biggl[1-\frac{x\ln x}{(x+1)\ln(x+1)}\biggr]\biggr\}\lim_{x\to0^+}\ln x\\
&=\lim_{x\to0^+} \biggl[1-\frac{x}{\ln(x+1)}\cdot\frac{\ln x}{x+1}\biggr]\lim_{x\to0^+}\ln x\\
&=-\infty.
\end{align*}
\par
The function $G(x)$ can also be rearranged as
\begin{equation}\label{f1(x)f2(x)-dfn}
  G(x)=\frac{x\ln x}{\ln(x+1)}[\ln(x+1)-\ln x]\triangleq f_1(x)f_2(x)
\end{equation}
for $x\in(0,\infty)$. It is not difficult to see that the function $f_2(x)$ is positive and decreasing on $(0,\infty)$. Straightforward differentiation produces
\begin{align*}
  f_1'(x)=\frac{\ln(x+1)+ [\ln(x+1)-{x}/{(x+1)}]\ln x}{[\ln(x+1)]^2} \triangleq\frac{g(x)}{[\ln(x+1)]^2}
\end{align*}
and $f_1(1)=0$. By virtue of the double inequality~\eqref{log-ineq-qi}, it follows that
\begin{align*}
  g(x)&>\frac{2x}{2+x}+\biggl[\frac{x(2+x)}{2(1+x)}-\frac{x}{x+1}\biggr]\ln x \\ &=\frac{2x}{2+x}+\frac{x^2}{2(1+x)}\ln x \\
  &=\frac{x^2}{2(1+x)}\biggl[\frac{4(1+x)}{x(2+x)}+\ln x\biggr]\\
  &\triangleq \frac{x^2}{2(1+x)}h(x)
\end{align*}
and
\begin{equation*}
  h'(x)=\frac{x^3-4x-8}{x^2 (x+2)^2}<0
\end{equation*}
for $x\in(0,1)$. As a result, the function $h(x)$ is decreasing on $(0,1)$ with $h(1)=\frac83$, so $h(x)>0$ on $(0,1)$. This means that the functions $g(x)$ and $f_1'(x)$ are positive on $(0,1)$. Hence, the function $f_1(x)$ is increasing and negative on $(0,1)$. In conclusion, the function $G(x)=f_1(x)f_2(x)$ is increasing and negative on $(0,1)$. The proof of Theorem~\ref{unit-ball-thm1} is complete.
\end{proof}

\begin{proof}[Second proof of a part of Theorem~\ref{unit-ball-thm1}]
By \eqref{f1(x)f2(x)-dfn}, it is clear that the function
$$
f_2(x)=\int_x^{x+1}\frac1t\td t=\int_0^1\frac1{t+x}\td t
$$
is completely monotonic on $(0,\infty)$.
\par
On the other hand, we have
\begin{align*}
f_1(x)&=\frac{u(x)-u(0)}{v(x)-v(0)}
\end{align*}
and
$$
\frac{u'(x)}{v'(x)}=(x+1)(1+\ln x)
$$
is increasing on $(0,\infty)$, where
$$
u(x)=\begin{cases}x\ln x,&x\in(0,1]\\0,&x=0
\end{cases}
$$
and $v(x)=\ln(x+1)$ for $x\in[0,1]$. The monotonic form of L'H\^ospital's rule put forward in~\cite[Theorem~1.25]{anderson1} (or see \cite[p.~92, Lemma~1]{Gene-Jordan-Inequal.tex} and \cite[p.~10, Lemma~2.9]{refine-jordan-kober.tex-JIA}) reads that if $U$ and $V$ are continuous on $[a,b]$ and differentiable on $(a,b)$ such that $V'(x)\ne0$ and $\frac{U'(x)}{V'(x)}$ is increasing $($or decreasing$)$ on $(a,b)$, then the functions $\frac{U(x)-U(b)}{V(x)-V(b)}$ and
$\frac{U(x)-U(a)}{V(x)-V(a)}$ are also increasing $($or decreasing$)$ on $(a,b)$. Therefore, the function $f_1(x)$ is increasing and negative on $(0,1)$. In a word, the function $G(x)=f_1(x)f_2(x)$ is increasing and negative on $(0,1)$.
\end{proof}

\end{document}